\documentclass[a4paper,12pt]{article}

\usepackage{amsfonts}

\newcommand{\C}{\mathbb{C}}
\newcommand{\Q}{\mathbb{Q}}

\newtheorem{thm}{Theorem}
\newtheorem{defn}[thm]{Definition}
\newtheorem{prop}[thm]{Proposition}
\newtheorem{lem}[thm]{Lemma}
\newtheorem{rem}[thm]{Remark}
\newtheorem{ex}[thm]{Example}
\newtheorem{exs}[thm]{Examples}
\newcommand{\dd}{\;\mathrm{d}}
\newcommand{\e}{\mathrm{e}}
\newcommand{\binom}[2]{{#1\choose#2}} 
\newcommand{\pf}{\noindent {\bf PROOF.} \quad}
\newcommand{\qed}{$\Box$}

\newcommand{\re}{\mathop{\mathrm{Re}}} 
\newcommand{\ii}{\mathrm{i}} 

\title{Higher Mahler measures and zeta functions}
\author{N. Kurokawa, M. Lal\'in\footnote{Supported by University of Alberta Fac. Sci. Startup Grant N031000610 and NSERC Discovery Grant 355412-2008},\,  and H. Ochiai}

\begin{document}
\baselineskip=17pt
\maketitle

\noindent{\bf Abstract}: We consider a generalization of the Mahler measure of a multivariable polynomial $P$ as the integral of $\log^k|P|$ in the unit torus, as opposed to the classical definition with the integral of $\log|P|$. A zeta Mahler measure, involving the integral of $|P|^s$, is also considered. Specific examples are computed, yielding special values of zeta functions, Dirichlet $L$-functions, and polylogarithms.

\medskip 



\section{Introduction} 
The (logarithmic) Mahler measure of a non-zero Laurent polynomial $P \in \C[x_1^{\pm1}, \dots, x_n^{\pm1}]$  is defined by
\[
m(P) = \int_0^1 \dots \int_0^1 
\log \left|P\left(\e^{2\pi\ii \theta_1}, \dots, \e^{2\pi\ii \theta_n}\right)\right| 
\dd \theta_1\cdots \dd \theta_n.\]

In this work, we consider the following generalization:
\begin{defn} The $k$-higher Mahler measure of $P$ is defined by
\[
m_k(P)
: =\int_0^1 \dots \int_0^1 
\log^k \left|P\left(\e^{2\pi\ii \theta_1}, \dots, \e^{2\pi\ii \theta_n}\right)\right| 
\dd \theta_1\cdots \dd \theta_n.
\]
\end{defn}

In particular, notice that for $k=1$ we obtain the classical Mahler measure
\[
m_1(P)=m(P),\quad \mathrm{and} \quad m_0(P)=1.
\]

These terms are the coefficients in the Taylor expansion of Akatsuka's zeta Mahler measure
\[ 
Z(s,P) 
= \int_0^1 \dots \int_0^1 
\left|P\left(\e^{2\pi\ii \theta_1}, \dots, \e^{2\pi\ii \theta_n}\right)\right|^s 
\dd \theta_1\dots \dd \theta_n,
\]
that is, 
\[Z(s,P) = \sum_{k=0}^\infty \frac{m_k(P) s^k}{k!}.\]
Akatsuka \cite{A} computed 
the zeta Mahler measure
$Z(s,x-c)$ for a constant $c$.


A natural generalization for the $k$-higher Mahler measure is the multiple higher Mahler measure for more than one polynomial.

\begin{defn} Let $P_1, \dots, P_l  \in \C[x_1^\pm, \dots, x_r^\pm]$ be non-zero Laurent polynomials. Their multiple higher Mahler measure is defined by
\begin{eqnarray*}
&& m(P_1,\dots,P_l) \\
&:=& \int_0^1 \cdots \int_0^1 
\left( \log\left|P_1(\e^{2\pi \ii \theta_1}, \cdots, \e^{2\pi \ii \theta_r}) \right| \right)
\cdots
\left( \log\left|P_l(\e^{2\pi \ii \theta_1}, \cdots, \e^{2\pi \ii \theta_r}) \right| \right)
\dd\theta_1 \cdots \dd\theta_r.
\end{eqnarray*}
\end{defn}

This construction yields the higher Mahler measures of one polynomial as a special case:
\[
m_k(P) = m(\underbrace{P,\ldots,P}_k).
\]

Moreover, the above definition implies that
\[
m(P_1)\cdots m(P_l) = m(P_1,\dots,P_l)
\]
when the variables of $P_j$'s in the right-hand side are algebraically independent.
This identity leads us to speculate about a product structure for the logarithmic Mahler measure. This would be a novel property, since the logarithmic Mahler measure is known to be additive, but no multiplicative structure is known.  

This definition has a natural counterpart in the world of zeta Mahler measures, namely, the {\em higher zeta Mahler measure} defined by 
\begin{eqnarray*}
&& Z(s_1,\dots,s_l; P_1,\dots,P_l) \\
&=& \int_0^1 \cdots \int_0^1 
\left|P_1(\e^{2\pi \ii \theta_1}, \cdots, \e^{2\pi \ii \theta_r}) \right|^{s_1}
\cdots
\left|P_l(\e^{2\pi \ii \theta_1}, \cdots, \e^{2\pi \ii \theta_r}) \right|^{s_l}
\dd\theta_1 \cdots \dd\theta_r,
\end{eqnarray*}

Its Taylor coefficients are related to the multiple higher Mahler measure:
\[
\frac{\partial^l}{\partial s_1 \cdots \partial s_l}
Z(0,\dots,0; P_1,\dots,P_l) = m(P_1,\dots,P_l).
\]

In this work, we compute the simplest examples of these heights and explore their basic properties. In section \ref{sec:high-onevar} we consider the case of higher Mahler measure for one-variable polynomials. More precisely, we consider linear polynomials in one variable. In particular, we obtain
\begin{eqnarray*}
m_2(x-1) &=& \frac{\pi^2}{12}, \\
m_3(x-1) &=& - \frac{3\zeta(3)}{2}, \\
m_4(x-1) &=& \frac{19 \pi^4}{240}, \\
m\left(1-x,1-\e^{2\pi\ii \alpha}x\right) & =& \frac{\pi^2}{2} \left( \alpha^2-\alpha+\frac{1}{6}\right), \qquad 0 \leq \alpha \leq 1. \\
\end{eqnarray*}
In section \ref{sec:m2}, we consider two examples of two-variable Mahler measure and we compute $m_2$. Sections \ref{sec:z}  and \ref{sec:zmulti} deal with examples of zeta Mahler measures of linear polynomials and their applications to the computation of higher Mahler measure, recovering the results from section \ref{sec:high-onevar} and giving an insight into them. Finally, we explore harder examples of zeta and higher Mahler measures in Section \ref{6}. For example,
\begin{eqnarray*}
m_2(x+y+2) & = &\frac{\zeta(2)}{2},\\
m_3(x+y+2) & = & \frac{9}{2} \log 2 \zeta(2)-\frac{15}{4} \zeta(3),\\
Z(s,x+x^{-1}+y+y^{-1}+c) &=& c^s _3F_2\left(\left.
\begin{array}{c} -\frac{s}{2}, \frac{1-s}{2}, \frac12 \\ 1,1 \end{array}
\right| \frac{16}{c^2} \right),\qquad c>4.\\
\end{eqnarray*}
\section{Higher Mahler measure of one-variable polynomials} \label{sec:high-onevar}
\subsection{The case of $1-x$}
Our first example is given by the simplest possible polynomial, namely $P=1-x$.
\begin{thm} \label{thm}
\[ 
m_k(1-x) 
=   \sum_{b_1+ \cdots + b_h = k,\, b_i \geq 2 } \frac{(-1)^kk!}{2^{2h}} \zeta(b_1, \dots, b_h),
\]
\end{thm}
where $\zeta(b_1, \dots, b_h)$ denotes a multizeta value, i.e.,
\[
\zeta(b_1, \dots, b_h) = \sum_{l_1< \dots <l_h}\frac{1}{l_1^{b_1}\dots l_h^{b_h}}.
\]
The right-hand side of Theorem~\ref{thm} 
can be re-written in terms of classical zeta values by using the following result.
\begin{prop} \label{prop} 
\begin{eqnarray*}
&&\sum_{\sigma \in S_h}\zeta(b_{\sigma(1)}, \dots, b_{\sigma(h)})
\\ 
&=& 
\sum_{ e_1+ \cdots + e_l =h } (-1)^{h-l} \prod_{s=1}^l (e_s-1)!  \sum  \zeta\left(\sum_{k \in \pi_1}b_k\right) \dots \zeta\left(\sum_{k \in \pi_l}b_k\right).
\end{eqnarray*}
where the sum in the right is taken over all the possible unordered partitions  of the set $\{1, \dots, h\}$ into $l$ subsets $\pi_1, \dots, \pi_l$ with $e_1, \dots, e_l$ elements respectively. 
\end{prop}

\pf (of Theorem~\ref{thm}) First observe that $x$ varies in the unit circle. Therefore, we can choose the principal branch for the logarithm. We proceed to write the function in terms of integrals of rational functions. We have
\begin{eqnarray*}
&& 
\log^k |1-x| = (\re \log (1-x))^k = \left(\frac{1}{2} (\log (1-x) + \log (1-x^{-1}))\right)^k 
\\ 
&=& 
\frac{1}{2^k} \left ( \int_0^1 \frac{\dd t}{t-x^{-1}} +\int_0^1 \frac{\dd t}{t-x} \right ) ^k 
= \frac{1}{2^k} \sum_{j=0}^k \binom{k}{j}
\left( \int_0^1 \frac{\dd t}{t-x^{-1}}\right)^j \left(  \int_0^1 \frac{\dd t}{t-x} \right) ^{k-j}.
\end{eqnarray*}
Now observe that 
\begin{eqnarray*} 
&& \left( \int_0^1 \frac{\dd t}{t-x^{-1}}\right)^j \left(  \int_0^1 \frac{\dd t}{t-x} \right ) ^{k-j} \\
&=& 
j! (k-j)! 
\int_0^1 \underbrace{\frac{\dd t}{t-x^{-1}} \circ \dots \circ \frac{\dd t}{t-x^{-1}}}_j 
\int_0^1 \underbrace{\frac{\dd t}{t-x} \circ \dots \circ \frac{\dd t}{t-x}}_{k-j}.
\end{eqnarray*}
We have just used the iterated integral notation of hyperlogarithms.

Combining the previous equalities gives 
\begin{eqnarray*}
&& m_k(1-x) = \frac{1}{2\pi \ii} \int_{|x|=1} \log^k |1-x| \frac{\dd x}{x} \\
&=& \frac{ k!}{2^k} \sum_{j=0}^k \frac{1}{2\pi \ii}  \int_{|x|=1} 
\int_0^1 \underbrace{\frac{\dd t}{t-x^{-1}} \circ \dots \circ \frac{\dd t}{t-x^{-1}}}_j 
\int_0^1 \underbrace{\frac{\dd t}{t-x} \circ \dots \circ \frac{\dd t}{t-x}}_{k-j} \frac{\dd x}{x}.
\end{eqnarray*}
If we now set $s=xt$ in the first $j$-fold integral and  $s=\frac{t}{x}$ in the second $(k-j)$-fold integral, the above becomes  
\[\frac{k!}{2^k} \sum_{j=0}^k  \frac{1}{2\pi \ii} \int_{|x|=1} 
\int_0^x \frac{\dd s}{s-1} \circ \dots \circ \frac{\dd s}{s-1} 
\int_0^{x^{-1}} \frac{\dd s}{s-1} \circ \dots \circ \frac{\dd s}{s-1} \frac{\dd x}{x}.
\]
We proceed to compute the integrals in terms of multiple polylogarithms:
\begin{eqnarray*}
m_k(1-x)&=& \frac{(-1)^kk!}{2^k} \sum_{j=0}^k \frac{1}{2\pi \ii}  \int_{|x|=1} \left( \sum_{0<l_1<\dots<l_{j}<\infty \,  0<m_1<\dots<m_{k-j} <\infty } \frac{x^{l_j-m_{k-j}}}{l_1\dots l_{j} m_1 \dots m_{k-j}} \right )  \frac{\dd x}{x}
\\
& =&  \frac{(-1)^kk!}{2^k} \sum_{j=1}^{k-1} \sum_{0<l_1<\dots<l_{j-1}<u <\infty,\,  0<m_1<\dots<m_{k-j-1}<u <\infty} \frac{1}{l_1\dots l_{j-1} m_1 \dots m_{k-j-1} u^2}.
\end{eqnarray*}

Now we need to analyze each term of the form
\begin{equation} \label{polysum}
\sum_{0<l_1<\dots<l_{j-1}<u <\infty,\,  0<m_1<\dots<m_{k-j-1}<u<\infty} \frac{1}{l_1\dots l_{j-1} m_1 \dots m_{k-j-1} u^2}.
\end{equation}

For an $h$-tuple $a_1, \dots, a_h$ such that $a_1+\dots +a_h = k-2h$, we set \[d_{a_1,\dots, a_h} = \sum_{e_1+\dots+e_h=j-h} \binom{a_1}{e_1} \dots \binom{a_h}{e_h} = \binom{a_1 + \dots + a_h}{e_1 + \dots +e_h} = \binom{k-2h}{j-h}.\]

Then the term (\ref{polysum}) is equal to
\[ \sum_{h=1}^{\min\{j-1,k-j-1\}} d_{a_1,\dots, a_h} \zeta( \{1\}_{a_1},2,\dots,\{1\}_{a_h},2).\]

Note that each term $\zeta( \{1\}_{a_1},2,\dots,\{1\}_{a_h},2)$ comes from choosing $h-1$ of the $l$'s and $h-1$ of the $m$'s and making them equal in pairs. Once this process has been done, one can choose the way the other $l$'s and $m$'s are ordered. All these choices give rise to the coefficients $d_{a_1, \dots, a_h}$.

The total sum is given by
\[m_k(1-x)= \sum_{h=1}^{k-1} c_{a_1,\dots, a_h} \zeta( \{1\}_{a_1},2,\dots,\{1\}_{a_h},2), \]
where
\[c_{a_1,\dots, a_h} = \frac{(-1)^kk!}{2^k} \sum_{j=1}^{k-1} \binom{k-2h}{j-h} = \frac{(-1)^kk!}{2^k}2^{k-2h} = \frac{(-1)^kk!}{2^{2h}}.\]

On the other hand,
\[\zeta( \{1\}_{a_1},2,\dots,\{1\}_{a_h},2) = \zeta(a_h+2, \dots ,a_1 +2).\]

To see this well-known fact, observe that the term in the left is
\[(-1)^{k-h}\int_0^1 \underbrace{\frac{\dd t}{t-1} \circ\dots \circ\frac{\dd t}{t-1}}_{a_1+1}\circ \frac{\dd t}{t}\circ \dots \circ\underbrace{\frac{\dd t}{t-1} \circ\dots \circ\frac{\dd t}{t-1}}_{a_h+1}\circ \frac{\dd t}{t}.\]
Making the change $t \rightarrow 1-t$ gives 
\[ = (-1)^{k-h}(-1)^k\int_0^1 \frac{\dd t}{t-1} \circ\underbrace{\frac{\dd t}{t} \circ\dots \circ\frac{\dd t}{t}}_{a_h+1}\circ \dots \circ \frac{\dd t}{t-1} \circ \underbrace{\frac{\dd t}{t} \circ\dots \circ\frac{\dd t}{t}}_{a_1+1}, \]
which corresponds to the term in the right. Thus, the total sum is
\[ m_k(1-x)= \sum_{b_1+ \dots + b_h = k,\, b_i \geq 2 } \frac{(-1)^kk!}{2^{2h}} \zeta(b_1, \dots, b_h).\]

\qed 

We show a proof of Proposition \ref{prop} for completeness.

\pf (Proposition \ref{prop}) We first show that we can write
\[ \sum_{\sigma \in S_h}\zeta(b_{\sigma(1)}, \dots, b_{\sigma(h)})\]
\[ = \sum_{ e_1+ \dots + e_l =h } r(e_1,\dots,e_l) \sum  \zeta\left(\sum_{k \in \pi_1}b_k\right) \dots \zeta\left(\sum_{k \in \pi_l}b_k\right).\]
where the function $r(e_1, \dots, e_l )$ satisfies some recurrence relationships. Here, as in the statement, the sum in the right is taken over all the possible unordered partitions  of the set $\{1, \dots, h\}$ into $l$ subsets $\pi_1, \dots, \pi_l$ with $e_1, \dots, e_l$ elements respectively.

Notice that $r$ is a function that is invariant under any permutation of its arguments.
We proceed by induction on $h$. It is clear that $r(1) =1$.  Also
\[ \zeta(a,b) + \zeta(b,a) = \zeta(a) \zeta(b) - \zeta(a+b), \]
from where $r(1,1) = 1$, $r(2) = -1$. 

Assume that the case of $h$ is settled. Now, we multiply everything by $\zeta(b_{h+1})$,
\[ \sum_{\sigma \in S_h}\zeta(b_{\sigma(1)}, \dots, b_{\sigma(h)}) \zeta(b_{h+1}) \]
\[ = \sum_{ e_1+ \dots + e_l =h } r(e_1,\dots,e_l) \sum  \zeta\left(\sum_{k \in \pi_1}b_k\right) \dots \zeta\left(\sum_{k \in \pi_l}b_k\right) \zeta(b_{h+1}).\]

Observe that 
\[ \sum_{\sigma \in S_h}\zeta(b_{\sigma(1)}, \dots, b_{\sigma(h)}) \zeta(b_{h+1})  = \sum_{\sigma \in S_{h+1}}\zeta(b_{\sigma(1)}, \dots, b_{\sigma(h+1)}) \]
\[ + \sum_{j=1}^h \sum_{\sigma \in S_h}\zeta(b_{\sigma(1)}, \dots, b_{\sigma(\breve{j})},\dots,b_{\sigma(h)}), \]
where $b_{\breve{j}} = b_j + b_{h+1}$. Hence,
\[ \sum_{\sigma \in S_{h+1}}\zeta(b_{\sigma(1)}, \dots, b_{\sigma(h+1)}) \]
\[ = \sum_{ e_1+ \dots + e_l =h } r(e_1,\dots,e_l) \sum  \zeta\left(\sum_{k \in \pi_1}b_k\right) \dots \zeta\left(\sum_{k \in \pi_l}b_k\right)\zeta(b_{h+1}) \]
\[ - \sum_{j=1}^h  \sum_{ e_1+ \dots + e_l =h } r(e_1,\dots,e_l) \zeta\left(\sum_{k \in \pi_1}b_k\right) \dots \zeta\left(b_{h+1} + \sum_{k \in \pi_f}b_k\right) \dots \zeta\left(\sum_{k \in \pi_l}b_k\right)\zeta(b_{h+1}).\]


From the above equation, we deduce the following identities:
\begin{eqnarray*}
r(e_1, \dots, e_f,1,e_{f+1},\dots, e_l )& = &r(e_1, \dots, e_f, e_{f+1},\dots , e_l ),\\
r(e_1, \dots, e_f + 1,\dots, e_l )& = &- e_f r(e_1, \dots, e_f, \dots  , e_l ).\\
\end{eqnarray*}

Now it is very easy to conclude that 
\begin{equation}
r(e_1, \dots, e_l )=(-1)^{h-l} \prod_{s=1}^l (e_s-1)!  .
\end{equation}
\qed

\begin{exs} \label{example}
Theorem~\ref{thm} enables us to compute $m_k(1-x)$. Here are the first few examples for $k=2,3,\dots,6$.

\begin{eqnarray*}
m_2(1-x) 
&=& \frac{ \zeta(2)}{2}, \\
m_3(1-x) 
&=& -6 \left( \frac{\zeta(3)}{4} \right) = - \frac{3 \zeta(3)}{2},  \\
m_4(1-x)
&=& 24 \left( \frac{\zeta(4)}{4} + \frac{\zeta(2,2)}{16} \right) \\
&=& 6 \zeta(4) + \frac{3 (\zeta(2)^2 - \zeta(4))}{4} = \frac{3 \zeta(2)^2 + 21\zeta(4)}{4}, \\
m_5(1-x)
&=& -120 \left( \frac{\zeta(5)}{4} + \frac{\zeta(2,3) + \zeta(3,2) }{16} \right) \\ 
&=& -30 \zeta(5) - \frac{15 (\zeta(2)\zeta(3) - \zeta(5))}{2} =- \frac{15 \zeta(2)\zeta(3) + 45\zeta(5)}{2}, \\
m_6(1-x)
&=& 720 \left( \frac{\zeta(6)}{4} + \frac{\zeta(3,3)}{16} +\frac{ \zeta(2,4)+ \zeta(4,2) }{16} + \frac{\zeta(2,2,2)}{64} \right) \\
&=& 180 \zeta(6) + \frac{45 (\zeta(3)^2 - \zeta(6))}{2} + 45( \zeta(2)\zeta(4) - \zeta(6)) \\
&& \qquad\qquad\qquad + \frac{45( 2 \zeta(6) - 3 \zeta(2)\zeta(4) + \zeta(2)^3) }{4 \cdot 6 } \\ 
&=& \frac{930 \zeta(6) + 180 \zeta(3)^2 + 315 \zeta(2)\zeta(4) + 15\zeta(2)^3}{8}.
\end{eqnarray*}
\end{exs}

\begin{rem}
Ohno and Zagier \cite{OZ} prove a 
result that generalizes Proposition~\ref{prop}.
Following their notation from (Theorem 1, \cite{OZ}), and setting $y=0$, $z=\frac{x^2}{4}$, (so that $s=n$) we have 
\[ \sum_{k=2}^\infty\sum_{b_1+ \dots + b_h = k,\, b_i \geq 2 } \frac{1}{2^{2h}} \zeta(b_1, \dots, b_h)x^k= \exp \left( \sum_{t=2}^\infty \frac{\zeta(t)}{t} x^t\left(1-\frac{1}{2^{t-1}}\right)\right). \]
\end{rem}
This identity also explains the relationship between the result in the statement of Theorem \ref{thm} and the result that is re obtained in Section \ref{4.2}.


\subsection{Higher Mahler measure for several linear polynomials}

As before, the simplest case to consider involves linear polynomials in one variable.
\begin{thm}
For $0 \le \alpha \le 1$
\[
m(1-x,1-\e^{2\pi \ii \alpha}x ) = \frac{\pi^2}{2}\left(\alpha^2-\alpha+\frac16 \right).
\]
\end{thm}
In particular, one obtains the following examples:
\begin{exs}
\begin{eqnarray*}
m(1-x,1-x) &=& \frac{\pi^2}{12}, \\
m(1-x,1+x) &=& - \frac{\pi^2}{24}, \\
m(1-x,1\pm \ii x) &=& - \frac{\pi^2}{96}, \\
m(1-x, 1-\e^{2\pi \ii \alpha}x ) &=& 0 \Leftrightarrow \alpha=\frac{3\pm\sqrt{3}}{6}. 
\end{eqnarray*}
\end{exs}
\pf By definition, 
\begin{eqnarray*}
m(1-x,1-\e^{2\pi \ii \alpha}x )
&=&
\int_0^1 \re \log(1-\e^{2\pi \ii \theta}) \cdot \re \log(1-\e^{2\pi \ii (\theta+\alpha)} ) \dd\theta \\
&=&
\int_0^1 
\left( - \sum_{k=1}^\infty \frac1k \cos 2 \pi k \theta \right)
\left(- \sum_{l=1}^\infty\frac1l \cos 2\pi (\theta+\alpha) \right) \dd\theta \\
&=&
\sum_{k,l \ge 1} \frac1{kl} \int_0^1 \cos (2\pi k \theta) \cos(2\pi l (\theta+\alpha)) \dd\theta.
\end{eqnarray*}
On the other hand,
\[
\int_0^1 \cos (2\pi k \theta) \cos(2\pi l (\theta+\alpha)) \dd\theta
= \left\{
\begin{array}{ll}
\frac12 \cos(2 \pi k \alpha) & \mbox{if } l = k, \\\\
0 & \mbox{otherwise}.
\end{array}\right.
\]
By putting everything together we conclude that
\begin{eqnarray*}
m(1-x,1-\e^{2\pi \ii \alpha}x )
&=& \frac12 \sum_{k=1}^\infty \frac{\cos(2\pi k \alpha)}{k^2} =\frac{\pi^2}{2}\left(\alpha^2-\alpha+\frac16\right).
\end{eqnarray*}
\qed

\begin{rem} \label{jensen}
The same calculation shows that
\[m(1-\alpha x, 1-\beta x ) =  \left\{
\begin{array}{ll}
\frac12 \re \mbox{\rm Li}_2\left(\alpha \bar{\beta}\right) &
\mbox{if } \left| \alpha \right|,  \left| \beta \right| \le 1, \\ \\
\frac12  \re \mbox{\rm Li}_2\left(\frac{\alpha \beta}{|\alpha|^2}\right) &
\mbox{if } \left| \alpha \right| \ge 1, \left| \beta \right| \le 1, \\ \\
\frac12  \re \mbox{\rm Li}_2\left( \frac{\alpha \bar{\beta}}{|\alpha\beta|^2}  \right)+\log|\alpha|\log|\beta| &
\mbox{if } \left| \alpha \right| , \left| \beta \right| \ge 1. \\
\end{array}\right.
\]

From this, one sees that for $P\in \C[x^\pm]$, $m_2(P)$ is a combination of dilogarithms and products of logarithms. 
In fact, for $P(x) = c x^s \prod_{j=1}^r (1-\alpha_j x )$,
we have 
\[m_2(P)=m(P,P)=(\log\left|c \right|)^2
+2(\log\left|c\right|) \sum_{j=1}^r \log^+|\alpha_j| +\sum_{j,k=1}^r m(1-\alpha_jx,1-\alpha_kx).\]
\end{rem}

The formula above plays an analogous role to Jensen's formula.

\begin{rem}
The previous computations may be extended to multiple higher Mahler measures involving more than two linear polynomials. For example,
\begin{eqnarray*}
m(1-x,1-\e^{2\pi i \alpha} x, 1-\e^{2 \pi i \beta} x) 
&=& - \frac14 \sum_{k,l\ge 1} \frac{\cos 2\pi((k+l)\beta-l \alpha)}{k l (k+l)} \\
&& 
- \frac14 \sum_{k,m\ge 1} \frac{\cos 2\pi((k+m)\alpha-m \beta)}{k m (k+m)} \\
&& - \frac14 \sum_{l,m\ge 1} \frac{\cos 2\pi(l \alpha+m\beta )}{l m(l+m)}. 
\end{eqnarray*}
\qed
\end{rem}

\section{Higher Mahler measure of two-variable polynomials} \label{sec:m2}

In this section we are going to consider examples of higher Mahler measures of polynomials in two variables. In particular, we will focus on the computation of $m_2$ using the formula from Remark \ref{jensen}, analogously to the way Jensen's formula for computing the classical Mahler measure of multivariable polynomials. 

The two polynomials that we consider were among the first examples of multivariable polynomials to be computed in terms of Mahler measure (by Smyth \cite{S}).

\subsection{$m_2(x+y+1)$}
\begin{thm}
\[m_2(x+y+1) = \frac{5 \pi^2}{54}.\]
\end{thm}

\pf 
We have, by definition,
\[ m_2(x+y+1) = \frac{1}{(2\pi\ii)^2}\int_{|y|=1} \int_{|x|=1} \log^2|x+y+1| \frac{\dd x}{x}\frac{\dd y}{y}.\]

We apply the result from Remark \ref{jensen} respect to the variable $y$,
\begin{eqnarray*}
m_2(x+y+1)&=& 
\frac{1}{2\pi\ii}\int_{|x|=1, |x+1|\leq 1} \frac{1}{2} \mbox{\rm Li}_2\left(|1+x|^2\right) \frac{\dd x}{x} 
\\
&+& \frac{1}{2\pi\ii}\int_{|x|=1, |x+1|\geq 1} \left( \frac{1}{2} \mbox{\rm Li}_2\left(\frac{1}{|1+x|^2}\right) + \log^2|1+x|\right ) \frac{\dd x}{x}.
\end{eqnarray*} 

Recalling the functional identity for the dilogarithm,
\[\mathrm{Li}_2(z) = -\mathrm{Li}_2\left(\frac{1}{z}\right) - \frac{1}{2} \log^2(-z)-\frac{\pi^2}{6}\]
for $z \not \in (0,1)$, we obtain
\begin{eqnarray*}
m_2(x+y+1)&=& \frac{1}{2\pi\ii}\int_{|x|=1, |x+1|\leq 1} \frac{1}{2} \mbox{\rm Li}_2\left(|1+x|^2\right) \frac{\dd x}{x}  \\
&+& \frac{1}{2\pi\ii}\int_{|x|=1, |x+1|\geq 1} \left( -\frac{1}{2} \re\mbox{\rm Li}_2\left(|1+x|^2\right) +\frac{\pi^2}{6}\right)\frac{\dd x}{x} \\
&=&
\frac{1}{2\pi\ii}\int_{|x|=1, |x+1|\leq 1} \mbox{\rm Li}_2\left(|1+x|^2\right) \frac{\dd x}{x}+\frac{ \pi^2}{9}
\\
&=& \frac{1}{2\pi}\int_{\frac{2\pi}{3}}^{\frac{4\pi}{3}}  \mbox{\rm Li}_2\left(4 \cos^2 \left(\frac{\theta}{2}\right)\right) 
\dd \theta + \frac{\pi^2}{9}.
\end{eqnarray*}

Notice that
\[\int \cos^{2n} \theta \dd \theta = \frac{\tan \theta}{2} \binom{2n-1}{n-1}\sum_{l=1}^{n} \frac{1} {2^{2n-2l}(2l-1)\binom{2l-2}{l-1}}\cos^{2l} \theta \]
\[ + \frac{1}{2^{2n-1}}\binom{2n-1}{n-1} \theta.\]

In particular,
\[\int_{\frac{\pi}{3}}^{\frac{2\pi}{3}} \cos^{2n} \theta \dd \theta = -\frac{\sqrt{3}}{2^{2n}} \binom{2n-1}{n-1}\sum_{l=0}^{n-1} \frac{1} {(2l+1)\binom{2l}{l}}  + \frac{1}{2^{2n-1}}\binom{2n-1}{n-1} \frac{\pi}{3}.\]

Now we use the identity for the sum of the inverses of Catalan numbers,
\[\frac{2 \pi \sqrt{3}}{9} =\sum_{l=0}^{\infty} \frac{1} {(2l+1)\binom{2l}{l}},\]
in order to get
\[\int_{\frac{\pi}{3}}^{\frac{2\pi}{3}} \cos^{2n} \theta \dd \theta =\frac{\sqrt{3}}{2^{2n}} \binom{2n-1}{n-1}\sum_{l=n}^\infty \frac{1} {(2l+1)\binom{2l}{l}}.\]


Note that
\[\frac{l! l!}{(2l+1)!} = B(l+1, l+1) =\int_0^1 s^l(1-s)^l \dd s.\]
Thus the above sum may be written as
\[\sum_{l=n}^\infty \int_0^1 s^l(1-s)^l \dd s = \int_0^1 \frac{s^n(1-s)^n}{1-s(1-s)} \dd s.\]

Putting everything together yields
\[\frac{1}{2\pi}\int_{\frac{2\pi}{3}}^{\frac{4\pi}{3}}  \mbox{\rm Li}_2\left(4 \cos^2 \left(\frac{\theta}{2}\right)\right) 
\dd \theta +\frac{\pi^2}{9}=\frac{\sqrt{3}}{\pi} \sum_{n=1}^\infty \frac{1}{n^2}\binom{2n-1}{n-1}\int_0^1 \frac{s^n(1-s)^n}{1-s(1-s)} \dd s +\frac{\pi^2}{9},\]
\begin{equation} \label{eq:1+x+yproof}
= \frac{\sqrt{3}}{{2\pi}} \int_0^1 \sum_{n=1}^\infty \frac{1}{n^2}\binom{2n}{n}\frac{s^n(1-s)^n}{1-s(1-s)} \dd s +\frac{\pi^2}{9}.
\end{equation}

At this point, we need the following
\begin{lem} \label{Lemma:Hiroyuki} For $|t|\leq \frac{1}{4}$, we have
 \begin{equation} \label{eq:Hiroyuki}
\sum_{k=1}^\infty \binom{2k}{k} \frac{t^k}{k^2}
= 2 \mathrm{Li}_2\left(\frac{1-\sqrt{1-4t}}{2}\right) - \left(\log\left(\frac{1+\sqrt{1-4t}}{2}\right)\right)^2.
\end{equation}

\end{lem}
\pf (of Lemma). We start from the series
\[
\sum_{k=1}^\infty \binom{2k}{k} t^k
= -1 + \sum_{k=0}^\infty \binom{-\frac12}{k} (-4t)^k
= -1 + \frac{1}{\sqrt{1-4t}}.
\]
convergent for $|t|\leq \frac{1}{4}$.

By integration,
we have
\[
\sum_{k=1}^\infty \binom{2k}{k} \frac{t^k}{k}
= -2 \log\left(1+\sqrt{1-4t}\right) + 2 \log 2.
\]
By integration again, we obtain the result.\qed

Now, if we set $t=s(1-s)$, we obtain $\frac{1-\sqrt{1-4t}}{2}=s$. Then the quantity $(\ref{eq:1+x+yproof})$ becomes
\[= \frac{\sqrt{3}}{{2\pi}} \int_0^1 (2 \mathrm{Li}_2(s) - \log^2 (1-s))\frac{\dd s}{1-s(1-s)} +\frac{\pi^2}{9}.\]
\[= -\frac{\sqrt{3}}{{\pi}} \int_{0\leq s_1\leq s_2 \leq s \leq 1} \frac{\dd s_1}{s_1-1} \frac{\dd s_2}{s_2} \frac{\dd s}{1-s+s^2} 
-\frac{\sqrt{3}}{{\pi}} \int_{0\leq s_1\leq s_2 \leq s \leq 1} \frac{\dd s_1}{s_1-1} \frac{\dd s_2}{s_2-1}\frac{\dd s}{1-s+s^2}  +\frac{\pi^2}{9}.\]
But 
\[\frac{1}{1-s+s^2} = \frac{1}{\sqrt{3} \ii} \left(\frac{1}{s-\omega} - \frac{1}{s-\bar{\omega}}\right),\]
where $\omega= \frac{1+\sqrt{3} \ii}{2}$. Thus, the above equals
\[ \frac{\ii}{{\pi}} \int_{0\leq s_1\leq s_2 \leq s \leq 1} \frac{\dd s_1}{s_1-1} \frac{\dd s_2}{s_2}  \left(\frac{1}{s-\omega} - \frac{1}{s-\bar{\omega}}\right) \dd s\]
\[+\frac{\ii}{{\pi}} \int_{0\leq s_1\leq s_2 \leq s \leq 1} \frac{\dd s_1}{s_1-1} \frac{\dd s_2}{s_2-1} \left(\frac{1}{s-\omega} - \frac{1}{s-\bar{\omega}}\right) \dd s +\frac{\pi^2}{9}.\]
\[= \frac{\ii}{{\pi}}(\mathrm{Li}_{2,1}(\omega, \bar{\omega}) - \mathrm{Li}_{2,1}(\bar{\omega}, \omega) -\mathrm{Li}_{1,1,1}(1,\omega,\bar{\omega}) + \mathrm{Li}_{1,1,1}(1,\bar{\omega}, \omega)) +\frac{\pi^2}{9}.\]
where we have written the result in terms of polylogarithms.

Now
\[\mathrm{Li}_{1,1,1}(1,\bar{\omega},\omega) -\mathrm{Li}_{1,1,1}(1,\omega,\bar{\omega})= \frac{5 \ii \pi^3}{81},\]
and
\[\mathrm{Li}_{2,1}(\bar{\omega},\omega) - \mathrm{Li}_{2,1}(\omega,\bar{\omega}) =\frac{7\ii \pi^3}{162}.\]
(see for example \cite{B}), and note that
 \[\frac{7 \pi^2}{162} -\frac{5 \pi^2}{81} +\frac{\pi^2}{9} =\frac{5\pi^2}{54}.\]

The result should be compared to Smyth's formula 
\[m(x+y+1) = \frac{3 \sqrt{3}}{4 \pi} L(\chi_{-3},2) = L'(\chi_{-3}, -1).\]

\subsection{$m_2(1+x+y(1-x))$}

\begin{thm}
 \[ m_2(1+x+y(1-x))= \frac{4 \ii}{\pi} (\mathrm{Li}_{2,1}(-\ii,-\ii) - \mathrm{Li}_{2,1}(\ii,\ii))+\frac{6\ii}{\pi}( -\mathrm{Li}_{2,1}(-\ii,\ii) + \mathrm{Li}_{2,1}(\ii,-\ii) )\]
\[+\frac{\ii}{\pi}(-\mathrm{Li}_{2,1}(1,\ii) +\mathrm{Li}_{2,1}(1,-\ii)  )-\frac{7\zeta(2)}{16}+\frac{\log 2}{\pi} L(\chi_{-4},2). \]
\end{thm}

\pf In order to apply the formula from Remark \ref{jensen} (for the variable $y$) we need to have a rational function that is monic in $y$. Therefore, we divide by the factor $(1+x)$:
\begin{equation}\label{eq3}
m_2(1-x+y(1+x))= m_2\left( \left(\frac{1-x}{1+x}\right)+y\right) + 2m\left( \left(\frac{1-x}{1+x}\right)+y, 1+x\right)+m_2(1+x).
\end{equation}

For the first term, we have
\[m_2\left( \left(\frac{1-x}{1+x}\right)+y\right)= \frac{1}{(2\pi \ii)^2} \int_{|y|=1}\int_{|x|=1} \log^2\left|\left(\frac{1-x}{1+x}\right)+y\right|\frac{\dd x}{x}\frac{\dd y}{y}.\]
By applying Remark \ref{jensen}, this equals
\[ \frac{1}{2\pi\ii} \int_{|x|=1, \, |1-x|\leq |1+x|}\frac{1}{2} \mathrm{Li}_2\left(\left|\frac{1-x}{1+x}\right|^2 \right)\frac{\dd x}{x} + \frac{1}{2\pi\ii} \int_{|x|=1, \, |1-x|\geq |1+x|}\frac{1}{2} \mathrm{Li}_2\left(\left|\frac{1+x}{1-x}\right|^2 \right)\frac{\dd x}{x} \]
\[+ \frac{1}{2\pi\ii} \int_{|x|=1, \, |1-x|\geq |1+x|} \log^2\left|\frac{1-x}{1+x} \right|\frac{\dd x}{x} \]
\[=\frac{1}{2\pi\ii} \int_{|x|=1, \, |1-x|\leq |1+x|}\mathrm{Li}_2\left(\left|\frac{1-x}{1+x}\right|^2 \right)\frac{\dd x}{x} + \frac{1}{2\pi\ii} \int_{|x|=1, \, |1-x|\geq |1+x|} \log^2\left|\frac{1-x}{1+x} \right|\frac{\dd x}{x}.\]

For the second term in equation (\ref{eq3}) we obtain
\[m\left( \left(\frac{1-x}{1+x}\right)+y,1+x\right)=\frac{1}{(2\pi\ii)^2} \int_{|y|=1}\int_{|x|=1} \log\left|\left(\frac{1-x}{1+x}\right)+y\right| \log|1+x| \frac{\dd x}{x}\frac{\dd y}{y}.\]
By Jensen's formula respect to the variable $y$, this equals
\[\frac{1}{2\pi\ii} \int_{|x|=1} \log^+\left|\frac{1-x}{1+x}\right| \log|1+x| \frac{\dd x}{x} =\frac{1}{2\pi\ii} \int_{|x|=1, \, |1-x|\geq |1+x|} \log\left|\frac{1-x}{1+x}\right|\log|1+x| \frac{\dd x}{x}.\]
Then (\ref{eq3}) becomes
\begin{eqnarray} \label{eq4}
m_2(1-x+y(1+x))&=& \frac{1}{2\pi\ii} \int_{|x|=1, \, |1-x|\leq |1+x|}\mathrm{Li}_2\left(\left|\frac{1-x}{1+x}\right|^2 \right)\frac{\dd x}{x} \nonumber \\
&&+ \frac{1}{2\pi\ii} \int_{|x|=1, \, |1-x|\geq |1+x|} (\log^2\left|1-x\right|-\log^2|1+x|)\frac{\dd x}{x} \nonumber\\
&&+ \frac{\zeta(2)}{2}.
\end{eqnarray}
For the first term on the right-hand side,
\begin{eqnarray*}
&& \frac{1}{2\pi\ii} \int_{|x|=1, \, |1-x|\leq |1+x|}\mathrm{Li}_2\left(\left|\frac{1-x}{1+x}\right|^2 \right)\frac{\dd x}{x} \\
&=& \frac{2}{\pi} \int_{-\frac{\pi}{4}}^{\frac{\pi}{4}}\mathrm{Li}_2\left(\tan^2 \theta  \right)\dd \theta = \frac{4}{\pi} \int_{-\frac{\pi}{4}}^{\frac{\pi}{4}}(\mathrm{Li}_2\left(\tan \theta  \right)+\mathrm{Li}_2\left(-\tan \theta  \right)) \dd \theta.
\end{eqnarray*}
After the change of variables $y=\tan \theta$, this becomes
\begin{eqnarray*}
\frac{8}{\pi} \int_{0}^{1}(\mathrm{Li}_2\left(y \right)+\mathrm{Li}_2\left(-y \right)) \frac{\dd y}{y^2+1}
&=& \frac{4}{\pi} \int_{0}^{1}(\mathrm{Li}_2\left(y \right)+\mathrm{Li}_2\left(-y \right)) \left(\frac{1}{1+\ii y}+\frac{1}{1-\ii y}\right)\dd y\\
&=& \frac{4}{\pi} (\ii \mathrm{Li}_{2,1}(\ii,-\ii) +\ii \mathrm{Li}_{2,1}(-\ii,-\ii) -\ii \mathrm{Li}_{2,1}(-\ii,\ii)-\ii \mathrm{Li}_{2,1}(\ii,\ii)).
\end{eqnarray*}
For the second term in (\ref{eq4}), we have
\begin{eqnarray*}
&& \frac{1}{2\pi\ii} \int_{|x|=1, \, |1-x|\geq |1+x|} (\log^2\left|1-x\right|-\log^2|1+x|)\frac{\dd x}{x} \\
&=&\sum_{k,l \geq 1}\frac{1-(-1)^{k+l}}{kl}2\int_{\frac{1}{4}}^{\frac{3}{4}} \cos(2\pi k \theta) \cos(2 \pi l \theta) \dd \theta\\ 
&=& \sum_{k,l \geq 1} \frac{1-(-1)^{k+l}}{2\pi kl}\left(\frac{\ii^{k+l+1}(1-(-1)^{k+l})}{k+l}+\frac{\ii^{k-l+1}(1-(-1)^{k-l})}{k-l}\right)\\
&=& \frac{\ii}{\pi} \sum_{k, l \geq 1} \frac{(1-(-1)^{k+l})\ii^{k+l}}{ kl^2} - \frac{\ii}{\pi} \sum_{k, l \geq 1} \frac{(1-(-1)^{k+l})\ii^{k+l}}{ (k+l)l^2} \\
&&+  \frac{2\ii}{\pi} \sum_{k>l \geq 1} \frac{(1-(-1)^{k+l})\ii^{k-l}}{ (k-l)l^2} - \frac{2\ii}{\pi} \sum_{k> l \geq 1} \frac{(1-(-1)^{k+l})\ii^{k-l}}{ kl^2}\\
&=&\frac{\ii}{\pi}(\mathrm{Li}_1(\ii)\mathrm{Li}_2(\ii)- \mathrm{Li}_1(-\ii)\mathrm{Li}_2(-\ii) -\mathrm{Li}_{2,1}(1,\ii) +\mathrm{Li}_{2,1}(1,-\ii)  )\\
&& +\frac{2\ii}{\pi}(\zeta(2)(\mathrm{Li}_1(\ii)- \mathrm{Li}_1(-\ii)) -\mathrm{Li}_{2,1}(-\ii,\ii) + \mathrm{Li}_{2,1}(\ii,-\ii)  )\\
&=&\frac{\ii}{\pi}(-\ii\log 2 L(\chi_{-4},2) -\frac{\pi \ii}{16}\zeta(2) -\mathrm{Li}_{2,1}(1,\ii) +\mathrm{Li}_{2,1}(1,-\ii)  )\\
&&+\frac{2\ii}{\pi}(\zeta(2)\frac{\pi \ii}{2} -\mathrm{Li}_{2,1}(-\ii,\ii) + \mathrm{Li}_{2,1}(\ii,-\ii)  ).
\end{eqnarray*}
Putting everything together in (\ref{eq4}), we obtain the final result
\begin{eqnarray*}
&& m_2(1-x+y(1+x)) \\
&=& \frac{4 \ii}{\pi} (\mathrm{Li}_{2,1}(-\ii,-\ii) - \mathrm{Li}_{2,1}(\ii,\ii))+\frac{6\ii}{\pi}( -\mathrm{Li}_{2,1}(-\ii,\ii) + \mathrm{Li}_{2,1}(\ii,-\ii)  )\\
&&+\frac{\ii}{\pi}(-\mathrm{Li}_{2,1}(1,\ii) +\mathrm{Li}_{2,1}(1,-\ii)  )-\frac{7\zeta(2)}{16}+\frac{\log 2}{\pi} L(\chi_{-4},2). 
\end{eqnarray*}
\qed

The previous result should be compared to (see \cite{S})
\[m(1-x+y(1+x))=\frac{2}{\pi} L(\chi_{-4},2). \]

\section{Zeta Mahler measures} \label{sec:z}
In this section, we consider zeta Mahler measures. We compute some examples and apply them to the computation of higher Mahler measures. 

\subsection{$Z(s,x-1)$}
As usual, we start with the linear polynomial $x-1$.
\begin{thm}\label{Z(x-1)}
\begin{eqnarray*}
Z(s,x-1) 
&=&
\int_0^1 (2 \sin \pi \theta)^s \dd\theta \\
&=& \exp \left( \sum_{k=2}^\infty \frac{(-1)^k (1-2^{1-k}) \zeta(k)}{k} s^k \right)
\end{eqnarray*}
around $s=0$.
\end{thm}

This result is a particular case of a formula obtained by Akatsuka \cite{A}.

\pf 
First we show that
\[
Z(s,x-1) 
= \frac{\Gamma(s+1)}{\Gamma(\frac s2+1)^2}
= \frac{s!}{((\frac s2)!)^2}
= \binom{s}{s/2},
\]
where $s!=\Gamma(s+1)$. In fact,
\[Z(s,x-1) 
= 2^{s+1} \int_0^{1/2} (\sin \pi \theta)^s \dd\theta. 
\]
After the change of variables $t=\sin^2 \pi \theta$ this becomes
\[ \frac{2^{s}}{\pi} \int_0^1 t^{\frac{s-1}{2}} (1-t)^{-1/2} \dd t.\]
So, we have obtained the Beta function:
\begin{eqnarray*}
Z(s, x-1)&=& \frac{2^s}{\pi} B\left(\frac{s+1}{2},\frac12\right) = \frac{2^s}{\pi} \frac{\Gamma\left(\frac{s+1}{2}\right) \Gamma\left(\frac12\right)}{\Gamma\left(\frac{s}{2}+1\right)} = \frac{2^s}{\sqrt{\pi}} \frac{\Gamma\left(\frac{s+1}{2}\right)}{\Gamma\left(\frac{s}{2}+1\right)}.
\end{eqnarray*}
Hence, by using
\[
\Gamma\left(\frac{s+1}{2}\right) = \frac{\Gamma(s)}{\Gamma\left(\frac{s}{2}\right)} 2^{1-s} \pi^{\frac12}
= \frac{\Gamma(s+1)}{\Gamma\left(\frac{s}{2}+1\right)} 2^{-s} \pi^{\frac12},
\]
we conclude that
\begin{equation}
Z(s,x-1) = \frac{\Gamma(s+1)}{\Gamma\left(\frac{s}{2}+1\right)^2}.
\end{equation}

On the other hand, the product expression
\[
\Gamma(s+1)^{-1} = \e^{\gamma s} \prod_{n=1}^\infty \left( 1+\frac{s}{n} \right) \e^{-\frac{s}{n}}
\]
yields
\begin{eqnarray*}
Z(s,x-1)
&=&
\prod_{n=1}^\infty \frac{\left(1+\frac{s}{2n}\right)^2}{1+\frac{s}{n}} \\
&=&
\exp\left( \sum_{n=1}^\infty
\left\{ 2 \log\left(1+\frac{s}{2n}\right) - \log\left(1+\frac{s}{n}\right) \right\} \right) \\
&=&
\exp\left( \sum_{k=1}^\infty \frac{(-1)^{k-1}}{k} \sum_{n=1}^\infty
\left\{ 2 \left(\frac{1}{2n}\right)^k - \frac{1}{n^k} \right\} s^k \right) \\
&=&
\exp\left( \sum_{k=2}^\infty \frac{(-1)^{k-1}}{k} \zeta(k) 
(2^{1-k}-1) s^k \right) \\
&=&
\exp\left( \sum_{k=2}^\infty 
\frac{(-1)^{k}(1-2^{1-k}) \zeta(k)}{k}  s^k \right).
\end{eqnarray*}
\qed

An analogous idea for evaluating $Z(s,P)$ appears in \cite{RV1}.

\subsection{$m_k(x-1)$} \label{4.2}
We can now use the evaluation of $Z(s,x-1)$ to re obtain the formula for $m_k(x-1)$. From Theorem~\ref{Z(x-1)},
\begin{eqnarray*}
Z(s,x-1)
&=& \exp\left( \frac{\zeta(2)}{4} s^2 - \frac{\zeta(3)}{4} s^3 + \frac{7\zeta(4)}{32} s^4 + \cdots \right) \\
&=& 1+ \frac{\zeta(2)}{4} s^2 - \frac{\zeta(3)}{4} s^3 
+ \left(\frac{7\zeta(4)}{32}+\frac{\zeta(2)^2}{32}\right) s^4 + \cdots. 
\end{eqnarray*}
On the other hand, by construction,
\[
Z(s,x-1) = 
1+ m_1(x-1) s + \frac12 m_2(x-1) s^2 + \frac16 m_3(x-1) s^3 + \frac1{24} m_4(x-1) s^4 + \cdots.
\]

Putting both identities together, we recover the result from Theorem \ref{thm}. In particular,
\begin{eqnarray*}
m_1(x-1) &=& 0, \\
m_2(x-1) &=& \frac{\zeta(2)}{2} = \frac{\pi^2}{12}, \\
m_3(x-1) &=& - \frac{3\zeta(3)}{2}, \\
m_4(x-1) &=& \frac34(7\zeta(4) + \zeta(2)^2) = \frac{19 \pi^4}{240}, \cdots \\
\end{eqnarray*}

\section{A computation of higher zeta Mahler measure}\label{sec:zmulti}

We compute the simplest example of a higher zeta Mahler measure and apply it to multiple higher Mahler measures. 
\begin{thm}
\begin{itemize}
\item[(i)]
\begin{eqnarray*}
Z(s,t;x-1,x+1)
&=& \int_0^1 \left| 2 \sin \pi \theta\right|^s \left| 2 \cos \pi \theta \right|^t \dd\theta \\
&=& \frac{\Gamma(s+1)\Gamma(t+1)}
{\Gamma\left(\frac{s}{2}+1\right) \Gamma\left(\frac{t}{2}+1\right) \Gamma\left(\frac{s+t}{2}+1\right)} \\
&=& \frac{s! t!}{\left(\frac s2\right)! \left(\frac t2\right)! \left(\frac{s+t}{2}\right)!} \\
&=& \prod_{n=1}^\infty \frac{\left(1+\frac{s}{2n}\right)\left(1+\frac{t}{2n}\right)\left(1+\frac{s+t}{2n}\right)}
{\left(1+\frac{s}{n}\right)\left(1+\frac{t}{n}\right)}.
\end{eqnarray*}
\item[(ii)]
\begin{eqnarray*}
Z(s,t;x-1,x+1)
&=&
\exp\left( \sum_{k=2}^\infty \frac{(-1)^k}{k} \zeta(k) \left\{
(1-2^{-k}) (s^k+t^k) - 2^{-k}(s+t)^k \right\} \right) \\
&\in & \Q[\pi^2, \zeta(3), \zeta(5),\dots][[s,t]]
\end{eqnarray*}
around $s=t=0$.
\item[(iii)]
\begin{eqnarray*}
m(\underbrace{x-1,\dots,x-1}_k,\underbrace{x+1,\dots,x+1}_l) &=&
\int_0^1 (\log \left| 2 \sin \pi \theta \right|)^k 
(\log \left| 2 \cos \pi \theta \right| )^l \dd \theta
\end{eqnarray*}
belongs to $\Q[\pi^2, \zeta(3), \zeta(5), \zeta(7),\dots]$
for integers $k,l\ge0$.
\end{itemize}
\end{thm}
\pf
\begin{itemize}
\item[\em(i)] By definition, 
\begin{eqnarray*}
Z(s,t;x-1,x+1)
&=& 2^{s+t} \int_0^1 (\sin \pi \theta)^s \left| \cos \pi \theta \right|^t \dd\theta \\
&=& 2^{s+t+1} \int_0^{1/2}(\sin \pi \theta)^s (\cos \pi \theta)^t \dd\theta. \\
\end{eqnarray*}
By the change of variables $u = \sin^2 \pi \theta$,
\begin{eqnarray*}
Z(s,t;x-1,x+1)&=& \frac{2^{s+t}}{\pi} \int_0^1 u^{\frac{s-1}{2}} (1-u)^{\frac{t-1}{2}} \dd u\\
&=& \frac{2^{s+t}}{\pi} B\left(\frac{s+1}{2}, \frac{t+1}{2}\right) \\
&=& \frac{2^{s+t}}{\pi} \frac{\Gamma\left(\frac{s+1}{2}\right) \Gamma\left(\frac{t+1}{2}\right)}{\Gamma\left(\frac{s+t}{2}+1\right)}.
\end{eqnarray*}
We now use again the identity
\[
\Gamma\left(\frac{z+1}{2}\right) = 2^{-z} \pi^{\frac12} \frac{\Gamma\left(z+1\right)}{\Gamma\left(\frac{z}{2}+1\right)},
\]
to get
\begin{eqnarray*}
Z(s,t;x-1,x+1)
&=& \frac{\Gamma(s+1)\Gamma(t+1)}{\Gamma\left(\frac{s}{2}+1\right) \Gamma\left(\frac{t}{2}+1\right) \Gamma\left(\frac{s+t}{2}+1\right)} 
\\
&=&
\prod_{n=1}^\infty \frac{\left(1+\frac{s}{2n}\right)\left(1+\frac{t}{2n}\right)\left(1+\frac{s+t}{2n}\right)}{\left(1+\frac{s}{n}\right)\left(1+\frac{t}{n}\right)}.
\end{eqnarray*}

\item[\em(ii)]
The above expression yields
\begin{eqnarray*}
&& Z(s,t;x-1,x+1) \\
&=&
\exp\left(
\sum_{n=1}^\infty \left\{
\log\left(1+\frac{s}{2n}\right) + \log\left(1+\frac{t}{2n}\right) + \log\left(1+\frac{s+t}{2n}\right) \right.\right. \\
&& \hskip3cm \left.\left. - \log\left(1+\frac{s}{n}\right)-\log\left(1+\frac{t}{n}\right)
\right\}\right) \\
&=&
\exp\left(
\sum_{k=1}^\infty \frac{(-1)^{k-1}}{k} 
\sum_{n=1}^\infty \left\{
\left(\frac{s}{2n}\right)^k + \left(\frac{t}{2n}\right)^k + \left(\frac{s+t}{2n}\right)^k - \left(\frac{s}{n}\right)^k-\left(\frac{t}{n}\right)^k
\right\}\right) \\
&=&
\exp\left(
\sum_{k=2}^\infty \frac{(-1)^{k-1}}{k} \zeta(k)
\left\{ 2^{-k} s^k + 2^{-k} t^k + 2^{-k} (s+t)^k - s^k - t^k 
\right\}\right) \\
&=&
\exp\left(
\sum_{k=2}^\infty \frac{(-1)^{k}}{k} \zeta(k)
\left\{ (1-2^{-k}) s^k + (1-2^{-k}) t^k - 2^{-k} (s+t)^k 
\right\}\right). \\
\end{eqnarray*}
This power series belongs to $\Q[\pi^2, \zeta(3), \zeta(5), \zeta(7),\dots][[s,t]]$.
\item[\em(iii)]
From {\em(ii)}, we see that
\[
\frac{\partial^{k+l}}{\partial s^k \partial t^l}
Z(0,0;x-1,x+1) \in \Q[\pi^2, \zeta(3), \zeta(5), \zeta(7),\dots],
\]
which is simply
\[
m(\underbrace{x-1,\dots,x-1}_k,\underbrace{x+1,\dots,x+1}_l) =
\int_0^1 (\log\left| 2 \sin \pi \theta\right|)^k
(\log \left| 2 \cos \pi \theta \right|)^l d\theta.
\]
\qed
\end{itemize}

\begin{ex}
In order to compute examples, we compare the terms of lowest degrees in the two expressions of $Z(s,t;x-1,x+1)$. On the one hand, we have
\begin{eqnarray*}
&& Z(s,t;x-1,x+1) \\
&=& \exp\left(
\frac{\zeta(2)}{2}\left(
\frac34(s^2+t^2)-\frac14(s+t)^2 \right)
- \frac{\zeta(3)}{3}\left(
\frac78(s^3+t^3)-\frac18(s+t)^3 \right) + 
(\mathrm{degree} \ge 4)
\right) \\
&=& \exp\left(
\frac{\zeta(2)}{4}\left(
s^2+t^2-st
\right)
- \frac{\zeta(3)}{8}\left(
2s^3+2t^3-s^2t-st^2
\right) + 
(\mathrm{degree} \ge 4)
\right).
\end{eqnarray*}
On the other hand,
\begin{eqnarray*}
&& Z(s,t;x-1,x+1) \\
&=&
1+
\left(\frac12 m(x-1,x-1) s^2 + \frac12 m(x+1,x+1) t^2 + m(x+1,x-1) st \right) \\
&& + \left( \frac16 m(x-1,x-1,x-1) s^3 + \frac16 m(x+1,x+1,x+1) t^3 \right. \\
&& + \left. \frac12 m(x-1,x-1,x+1) s^2t  +  \frac12 m(x-1,x+1,x+1) st^2 \right) \\
&& + (\mathrm{degree} \ge 4).
\end{eqnarray*}
We obtain:
\begin{eqnarray*}
m(x-1,x+1) 
&=& \int_0^1 \log \left| 2 \sin \pi \theta \right| \log \left| 2 \cos \pi \theta \right| d \theta
= -\frac{\zeta(2)}{4} = - \frac{\pi^2}{24}, \\
m(x-1,x-1,x+1) 
&=&
\int_0^1 (\log \left| 2 \sin \pi \theta \right|)^2 \log \left| 2 \cos \pi \theta \right| d \theta
= 2\frac{\zeta(3)}{8} =  \frac{\zeta(3)}{4}, \\
m(x-1,x+1,x+1) 
&=&
\int_0^1 \log \left| 2 \sin \pi \theta \right| (\log \left| 2 \cos \pi \theta \right|)^2 d \theta
= 2\frac{\zeta(3)}{8} =  \frac{\zeta(3)}{4}.
\end{eqnarray*}
Note that the calculation
\[
Z(s,0;x-1,x+1) = Z(s,x-1) = \binom{s}{s/2}
\]
yields $m_k(x-1)$ again.

We also remark that we have another relation
\[
Z(s,s;x-1,x+1) = Z(s,x-1)=Z(s,x+1).
\]
\qed
\end{ex}

\section{Further examples} \label{6}

\subsection{The case $P=x+x^{-1}+y+y^{-1}+c$}
\begin{thm} \label{thm14}
For $c>4$,
\begin{eqnarray*}
 Z(s,x+x^{-1}+y+y^{-1}+c)
&=&  c^s \sum_{j=0}^\infty \binom{s}{2j} \frac{1}{c^{2j}} \binom{2j}{j}^2 \\
&=&c^s {}_3F_2\left(\left.
\begin{array}{c} -\frac{s}{2}, \frac{1-s}{2}, \frac12 \\ 1,1 \end{array}
\right| \frac{16}{c^2}
\right), \\
\end{eqnarray*}
where
the generalized hypergeometric series ${}_3F_2$ is defined by
\[
{}_3F_2\left(\left.\begin{array}{c} a_1,a_2,a_3 \\ b_1,b_2 \end{array} \right| z\right)
= \sum_{j=0}^\infty \frac{(a_1)_j (a_2)_j (a_3)_j}{(b_1)_j (b_2)_j j!} z^j,
\]
with the Pochhammer symbol defined by $(a)_j = a(a+1)\cdots (a+j-1)$.
\end{thm}
\pf
We first write $x+x^{-1}+y+y^{-1}+c = c\left( \frac{x+ x^{-1}+y+y^{-1}}{c}+1\right)$. Since $c \geq 4$, $ \frac{x+ x^{-1}+y+y^{-1}}{c}+1$ is a positive number in the unit torus. Hence, we may omit the absolute value in the computation of the zeta function. Therefore we may write
\begin{eqnarray*}
&&
 Z(s,x+x^{-1}+y+y^{-1}+c) \\
&=& \frac{1}{(2\pi\ii)^2}\int_{|y|=1} \int_{|x|=1} (x+x^{-1}+y+y^{-1}+c)^s \frac{\dd x}{x} \frac{\dd y}{y}\\
&=& \frac{c^s}{(2\pi\ii)^2}\int_{|y|=1} \int_{|x|=1} \left(1+\frac{x+x^{-1}+y+y^{-1}}{c}\right)^s \frac{\dd x}{x} \frac{\dd y}{y}\\
&=&c^s \sum_{k=0}^\infty \binom{s}{k} \frac{1}{(2\pi\ii)^2}\int_{|y|=1} \int_{|x|=1} \left(\frac{x+x^{-1}+y+y^{-1}}{c}\right)^k \frac{\dd x}{x} \frac{\dd y}{y}\\
&=&c^s \sum_{j=0}^\infty \binom{s}{2j} \frac{1}{c^{2j}} \binom{2j}{j}^2.\\
\end{eqnarray*}

The last equality is the result of the following observation. The number 
\[\frac{1}{(2\pi\ii)^2}\int_{|y|=1} \int_{|x|=1} \left(x+x^{-1}+y+y^{-1}\right)^k \frac{\dd x}{x} \frac{\dd y}{y}\]
is the constant coefficient of $\left(x+x^{-1}+y+y^{-1}\right)^k$. This idea was observed by Rodriguez-Villegas \cite{RV} who studied this specific example as part of the computation of the classical Mahler measure for this family of polynomials. 

The expression in terms of the generalized hypergeometric function 
is derived by
$\binom{s}{2j} (2j)!= 2^{2j} (-\frac{s}{2})_j (\frac{1-s}{2})_j$
and $(2j)! = 2^{2j} (\frac12)_j j!$.
Note that the series ${}_3F_2(z)$ converges in $\left| z \right|<1$,
which is compatible with the condition $c>4$ in the statement of the Theorem.
\qed

\subsection{Properties of zeta Mahler measures}
The proof of Theorem \ref{thm14} may also be achieved by combining the following elementary properties
of zeta Mahler measures:
\begin{lem}\label{property}
\begin{itemize}
\item[(i)]
For a positive constant $\lambda$,
we have $Z(s, \lambda P) = \lambda^s Z(s,P)$.
\item[(ii)]
Let $P \in \C[x_1^{\pm1},\dots,x_n^{\pm1}]$ be a Laurent polynomial such that it takes
non-negative real values in the unit torus.
Then we have the following series expansion on 
$\left| \lambda\right| \le 1/\max(P)$, where
$\max(P)$ is the maximum of $P$ on the unit torus:
\begin{eqnarray*}
Z(s,1+\lambda P) &=& \sum_{k=0}^\infty \binom{s}{k} Z(k, P) \lambda^k, \\
m(1+\lambda P) &=& \sum_{k=1}^\infty \frac{(-1)^{k-1}}{k} Z(k,P) \lambda^k. \\
\end{eqnarray*}
More generally,
\begin{eqnarray*}
m_j(1+\lambda P) &= &j!\sum_{0<k_1< \dots < k_j}  \frac{(-1)^{k_j-j} }{k_1 \dots k_j}Z(k_j,P)\lambda^{k_j}.
\end{eqnarray*}
\item[(iii)]
$Z(s,P) = Z(\frac{s}{2}, P \bar{P})$,
where we put $\bar{P} = \sum_{\alpha} \bar{a}_\alpha x^{-\alpha}$
for $P=\sum_{\alpha} a_\alpha x^\alpha$.
Note that $P \bar{P}$ is real-valued on the torus.
\end{itemize}
\end{lem}
Therefore, in principle, the knowledge of $m(1+\lambda P)$
yields enough information to determine $Z(s,1+\lambda P)$.

\pf 
{\em(i)} and {\em(iii)} are obvious.
For {\em (ii)}, 
we may use the Taylor expansions in $\lambda$;
\[(1+\lambda P)^s = \sum_{k=0}^\infty \binom{s}{k} \lambda^k P^k, \quad \log(1+\lambda P) = \sum_{k=1}^\infty \frac{(-1)^{k-1}}{k} \lambda^k P^k.\] 

In particular, we may write 
\[
Z(s,1+\lambda P) 
= \sum_{k=0}^\infty m_k(1+\lambda P) \frac{s^k}{k!} 
= \sum_{k=0}^\infty Z(k,P) \lambda^k \frac{s(s-1)\cdots(s-k+1)}{k!}.\]
In other words,
the coefficients with respect to the monomial basis are
the $k$-logarithmic Mahler measures $m_k(1+\lambda P)$,
while the coefficients with respect to the shifted monomial basis
are (the special values of) zeta Mahler measures $Z(k,P) \lambda^k$.

Combining these observations, we obtain the three equalities. \qed

\subsection{The case $P=x+y+c$}
Now we apply these ideas to $P=x+y+c$ with $c\geq 2$. 

\begin{thm}\label{x+y+c}
Let $c \ge 2$. Then 
\begin{itemize}
\item[(i)]
\begin{eqnarray*}
Z(s,x+y+c)
&=&c^s \sum_{j=0}^\infty \binom{s/2}{j}^2 \frac{1}{c^{2j}} \binom{2j}{j},
\end{eqnarray*}
\item[(ii)] 
\[
m_2(x+y+c) =\log^2 c +\frac{1}{2} \sum_{k=1}^\infty \binom{2k}{k}\frac{1}{k^2c^{2k}},
\]
\item[(iii)]
\[m_3(x+y+c)= \log^3 c + \frac{3}{2} \log c \sum_{k=1}^\infty \binom{2k}{k}\frac{1}{k^2 c^{2k}}-\frac{3}{2} \sum_{k=2}^\infty \binom{2k}{k}\frac{1}{k^2 c^{2k}} \sum_{j=1}^{k-1} \frac{1}{j}.\]
In particular, we obtain the special values
\item[(iv)]
 \[m_2(x+y+2)=\frac{\zeta(2)}{2},\]
\item[(v)]
\[m_3(x+y+2)= \frac{9}{2} \log 2 \zeta(2)-\frac{15}{4} \zeta(3).\]
\end{itemize}
\end{thm}
\pf 
{\em (i)}
In this case, the polynomial is not reciprocal, so we first need to consider $(x+y+c)(x^{-1}+y^{-1}+c)$. 
Then,
\begin{eqnarray*}
&& Z(s,x+y+c) \\
&=&Z(s/2, (x+y+c)(x^{-1}+y^{-1}+c))\\
&=& \frac{1}{(2\pi\ii)^2}\int_{|y|=1} \int_{|x|=1} ((x+y+c)(x^{-1}+y^{-1}+c))^{s/2} \frac{\dd x}{x} \frac{\dd y}{y}\\
&=& \frac{c^{s}}{(2\pi\ii)^2}\int_{|y|=1} \int_{|x|=1} \left(1+\frac{x+y}{c}\right)^{s/2}\left(1+\frac{x^{-1}+y^{-1}}{c}\right)^{s/2}\frac{\dd x}{x} \frac{\dd y}{y}\\
&=&c^s \sum_{j=0}^\infty\sum_{k=0}^\infty \binom{s/2}{j}  \binom{s/2}{k} \frac{1}{(2\pi\ii)^2}\int_{|y|=1} \int_{|x|=1} \left(\frac{x+y}{c}\right)^{j}\left(\frac{x^{-1}+y^{-1}}{c}\right)^{k} \frac{\dd x}{x} \frac{\dd y}{y}\\
&=&c^s \sum_{j=0}^\infty \binom{s/2}{j}^2 \frac{1}{c^{2j}} \binom{2j}{j}.   \\
\end{eqnarray*}
The last identity was obtained, as in the case of $x+x^{-1}+y+y^{-1}+c$, by computing the constant coefficient of the product of powers of polynomials in the integrand. 

Formulas {\em(ii)} and {\em(iii)} are consequence of {\em(i)} and Lemma \ref{property}. 

If we set $t=1/4$ in the equation of Lemma \ref{Lemma:Hiroyuki}, we obtain $\zeta(2) - 2 (\log 2)^2$. Combining this with {\em (ii)}, we get the result of {\em(iv)}.

For the last formula {\em(v)},
it is enough to prove the following identity:
\[\sum_{k=2}^\infty \binom{2k}{k}\frac{1}{k^2 4^{k}} \sum_{j=1}^{k-1} \frac{1}{j}= -\frac{4}{3}\log^3 2-2\zeta(2)\log 2 +\frac{5}{2}\zeta(3).\] 

We have
\[
\sum_{k=1}^\infty \binom{2k}{k} \frac{t^k}{k^2}
= 2 \mathrm{Li}_2\left(\frac{1-\sqrt{1-4t}}{2}\right) - \left(\log\left(\frac{1+\sqrt{1-4t}}{2}\right)\right)^2.
\]
Now we turn the left-hand side into a double series:
\begin{eqnarray*}
&&\sum_{k=2}^\infty \binom{2k}{k} \frac{t^k}{k^2} \sum_{j=0}^{k-2} x^j = \sum_{k=2}^\infty \binom{2k}{k} \frac{t^k}{k^2} \left(\frac {x^{k-1}-1}{x-1}\right) \\
&=& \frac{1}{x(x-1)}\left( 2\mathrm{Li}_2\left(\frac{1-\sqrt{1-4xt}}{2}\right) - \left(\log\left(\frac{1+\sqrt{1-4xt}}{2}\right)\right)^2\right)
\\
&&- \frac{1}{x-1}\left( 2\mathrm{Li}_2\left(\frac{1-\sqrt{1-4t}}{2}\right) - \left(\log\left(\frac{1+\sqrt{1-4t}}{2}\right)\right)^2\right).
\end{eqnarray*}
In particular, by evaluating at $t = \frac{1}{4}$, we obtain
\begin{eqnarray*}
&& \sum_{k=2}^\infty \binom{2k}{k} \frac{1}{k^24^k} \sum_{j=0}^{k-2} x^j \\
&=& \frac{1}{x(x-1)}\left(2\mathrm{Li}_2\left(\frac{1-\sqrt{1-x}}{2}\right) - \left(\log\left(\frac{1+\sqrt{1-x}}{2}\right)\right)^2\right)
\\
&& \hskip2cm - \frac{1}{x-1}\left( \zeta(2) -2 \log^2 2\right).
\end{eqnarray*}
Integrating from 0 to 1, we obtain the double series that we wish to evaluate:
\begin{eqnarray*}
I&:=&
\sum_{k=2}^\infty \binom{2k}{k} \frac{1}{k^24^k} \sum_{j=1}^{k-1} \frac{1}{j} \\
&=& \int_0^1 \left(\frac{1}{x(x-1)}\left(2\mathrm{Li}_2\left(\frac{1-\sqrt{1-x}}{2}\right) - \left(\log\left(\frac{1+\sqrt{1-x}}{2}\right)\right)^2\right)\right.
\\
&& \hskip2cm \left. - \frac{1}{x-1}\left( \zeta(2) -2 \log^2 2\right)\right) \dd x.
\end{eqnarray*}
We just need to perform the integration. For that, we consider the change of variables $y=\frac{1-\sqrt{1-x}}{2}$:
\begin{eqnarray*}
I&=& \int_0^\frac{1}{2} \left(2\mathrm{Li}_2\left(y\right) - \left(\log\left(1-y\right)\right)^2\right)\left(\frac{4}{2y-1}-\frac{1}{y-1}-\frac{1}{y} \right)\dd y
\\
&& - 4\left( \zeta(2) -2 \log^2 2\right)\int_0^\frac{1}{2} \frac{\dd y}{2y-1}.
\end{eqnarray*}
We write the expression in terms of iterated integrals, so that we can relate the result to multiple polylogarithms:
\begin{eqnarray*}
&& 2\mathrm{Li}_2\left(y\right) - \left(\log\left(1-y\right)\right)^2 \\
&=& -2 \int_{0\leq t_1\leq t_2\leq y} \frac{\dd t_1}{t_1-1} \frac{\dd t_2}{t_2} - 2\int_{0\leq t_1\leq t_2\leq y} \frac{\dd t_1}{t_1-1} \frac{\dd t_2}{t_2-1}.
\end{eqnarray*}
We have
\begin{eqnarray*}
I &=&  -2 \int_{0\leq t_1\leq t_2\leq y\leq \frac{1}{2}} \frac{\dd t_1}{t_1-1} \frac{\dd t_2}{t_2}\left(\frac{4}{2y-1}-\frac{1}{y-1}-\frac{1}{y} \right)\dd y \\
&& - 2\int_{0\leq t_1\leq t_2\leq y\leq \frac{1}{2}} \frac{\dd t_1}{t_1-1} \frac{\dd t_2}{t_2-1}\left(\frac{4}{2y-1}-\frac{1}{y-1}-\frac{1}{y} \right)\dd y\\
&& + \left(2 \int_{0\leq t_1\leq t_2\leq \frac{1}{2}} \frac{\dd t_1}{t_1-1} \frac{\dd t_2}{t_2} + 2\int_{0\leq t_1\leq t_2\leq  \frac{1}{2}} \frac{\dd t_1}{t_1-1} \frac{\dd t_2}{t_2-1} \right)\int_0^\frac{1}{2} \frac{4\dd y}{2y-1}.
\end{eqnarray*}
After some rearranging we get
\begin{eqnarray*}
I&= & 2 \int_{0\leq t_1\leq t_2\leq y\leq \frac{1}{2}} \frac{\dd t_1}{t_1-1} \frac{\dd t_2}{t_2}\left(\frac{1}{y-1}+\frac{1}{y} \right)\dd y \\
&& + 2\int_{0\leq t_1\leq t_2\leq y\leq \frac{1}{2}} \frac{\dd t_1}{t_1-1} \frac{\dd t_2}{t_2-1}\left(\frac{1}{y-1}+\frac{1}{y} \right)\dd y
\\
&&+ 8 \int_{0 \leq y , t_1\leq t_2\leq \frac{1}{2}} \frac{\dd y}{2y-1}\frac{\dd t_1}{t_1-1} \frac{\dd t_2}{t_2} + 8\int_{0\leq y, t_1\leq t_2\leq \frac{1}{2}} \frac{\dd y}{2y-1} \frac{\dd t_1}{t_1-1} \frac{\dd t_2}{t_2-1}.
 \end{eqnarray*}
We make the change of variables $s_i=2t_i$, $z=2y$. Then
\begin{eqnarray*}
I&=&  2 \int_{0\leq s_1\leq s_2\leq z\leq 1} \frac{\dd s_1}{s_1-2} \frac{\dd s_2}{s_2}\left(\frac{1}{z-2}+\frac{1}{z} \right)\dd z \\
&& + 2\int_{0\leq s_1\leq s_2\leq z\leq 1} \frac{\dd s_1}{s_1-2} \frac{\dd s_2}{s_2-2}\left(\frac{1}{z-2}+\frac{1}{z} \right)\dd z
\\
&& + 4 \int_{0 \leq z , s_1\leq s_2\leq 1} \frac{\dd z}{z-1}\frac{\dd s_1}{s_1-2} \frac{\dd s_2}{s_2} + 4\int_{0\leq z, s_1\leq s_2\leq 1} \frac{\dd z}{z-1} \frac{\dd s_1}{s_1-2} \frac{\dd s_2}{s_2-2}. 
\end{eqnarray*}
Now we make another change of variables $u_i=1-s_i$, $w=1-z$ to get
\begin{eqnarray*}
I&=&  -2 \int_{0\leq w\leq u_2\leq u_1\leq 1} \left(\frac{1}{w+1}+\frac{1}{w-1} \right)\dd w\frac{\dd u_2}{u_2-1} \frac{\dd u_1}{u_1+1} \\
&& - 2\int_{0\leq w\leq u_2\leq u_1\leq 1} \left(\frac{1}{w+1}+\frac{1}{w-1} \right)\dd w\frac{\dd u_2}{u_2+1}\frac{\dd u_1}{u_1+1}
\\
&&- 4 \int_{0 \leq u_2 \leq w, u_1 \leq 1} \frac{\dd u_2}{u_2-1}\frac{\dd u_1}{u_1+1} \frac{\dd w}{w} - 4\int_{0 \leq u_2 \leq w, u_1 \leq 1} \frac{\dd u_2}{u_2+1} \frac{\dd u_1}{u_1+1} \frac{\dd w}{w}.\\
\end{eqnarray*}
We may now express all the terms as hyperlogarithms, and then as multiple polylogarithms evaluated in $\pm 1$.
\begin{eqnarray*}
I&=& -2I_{1,1,1}(-1,1,-1,1)-2I_{1,1,1}(1,1,-1,1)- 2I_{1,1,1}(-1,-1,-1,1)-2I_{1,1,1}(1,-1,-1,1)\\
&&- 4 I_{1,2}(1,-1,1) - 4I_{2,1}(1,-1,1)-4 I_{1,2}(-1,-1,1) - 4I_{2,1}(-1,-1,1)\\
&=&2\mathrm{Li}_{1,1,1}(-1,-1,-1)+2\mathrm{Li}_{1,1,1}(1,-1,-1)+2\mathrm{Li}_{1,1,1}(1,1,-1)+2\mathrm{Li}_{1,1,1}(-1,1,-1)\\
&&-4\mathrm{Li}_{1,2}(-1,-1) - 4\mathrm{Li}_{2,1}(-1,-1)-4\mathrm{Li}_{1,2}(1,-1)-4\mathrm{Li}_{2,1}(1,-1).
\end{eqnarray*}
The terms involving multiple polylogarithms of length greater than 1 may be expressed as terms involving ordinary polylogarithms (of length 1). First, we reduce the multiple polylogarithms from length 3 to length 2 and 1 using the following identities:
\begin{eqnarray*}
\mathrm{Li}_{1,1,1}(-1,-1,-1)&=& \frac{1}{3}(\mathrm{Li}_1(-1)\mathrm{Li}_{1,1}(-1,-1)- \mathrm{Li}_{2,1}(1,-1) - \mathrm{Li}_{1,2}(-1,1)),\\
\mathrm{Li}_{1,1,1}(1,-1,-1)&=&\frac{1}{12} ( 6 \mathrm{Li}_1(-1)\mathrm{Li}_{1,1}(1,-1)-2\mathrm{Li}_{1}(-1) \mathrm{Li}_{1,1}(-1,-1)\\
&&\hskip-15mm
-\mathrm{Li}_{2,1}(1,-1)-\mathrm{Li}_{1,2}(-1,1) -6 \mathrm{Li}_{2,1}(-1,-1)-6\mathrm{Li}_{1,2}(1,1)),
\\
\mathrm{Li}_{1,1,1}(1,1,-1)&=&\frac{(\mathrm{Li}_1(-1))^3}{6}, \\
\mathrm{Li}_{1,1,1}(-1,1,-1) &=& \frac{1}{6}(2\mathrm{Li}_1(-1)\mathrm{Li}_{1,1}(-1,-1)+\mathrm{Li}_{2,1}(1,-1)+\mathrm{Li}_{1,2}(-1,1)).
\end{eqnarray*}
Incorporating these identities in the expression for $I$, we get
\begin{eqnarray*}
I&=&\frac{2}{3} \mathrm{Li}_1(-1)\mathrm{Li}_{1,1}(-1,-1)-\frac{2}{3} \mathrm{Li}_{2,1}(1,-1) -\frac{2}{3} \mathrm{Li}_{1,2}(-1,1)\\
&& +\mathrm{Li}_1(-1)\mathrm{Li}_{1,1}(1,-1)-\frac{1}{3}\mathrm{Li}_{1}(-1) \mathrm{Li}_{1,1}(-1,-1)-\frac{1}{6}\mathrm{Li}_{2,1}(1,-1)\\ &&-\frac{1}{6}\mathrm{Li}_{1,2}(-1,1) - \mathrm{Li}_{2,1}(-1,-1)-\mathrm{Li}_{1,2}(1,1)+\frac{1}{3}(\mathrm{Li}_1(-1))^3\\
&& +\frac{2}{3}\mathrm{Li}_1(-1)\mathrm{Li}_{1,1}(-1,-1)+\frac{1}{3}\mathrm{Li}_{2,1}(1,-1)+\frac{1}{3}\mathrm{Li}_{1,2}(-1,1)\\
&& -4\mathrm{Li}_{1,2}(-1,-1) - 4\mathrm{Li}_{2,1}(-1,-1)-4\mathrm{Li}_{1,2}(1,-1)-4\mathrm{Li}_{2,1}(1,-1) \\
& =& \mathrm{Li}_1(-1)\mathrm{Li}_{1,1}(-1,-1)-\frac{9}{2} \mathrm{Li}_{2,1}(1,-1) -\frac{1}{2} \mathrm{Li}_{1,2}(-1,1)\\
&&+\mathrm{Li}_1(-1)\mathrm{Li}_{1,1}(1,-1)- 5 \mathrm{Li}_{2,1}(-1,-1)-\mathrm{Li}_{1,2}(1,1)\\
&& +\frac{1}{3}(\mathrm{Li}_1(-1))^3-4\mathrm{Li}_{1,2}(-1,-1) -4\mathrm{Li}_{1,2}(1,-1).
\end{eqnarray*}
Now we consider identities of multiple polylogarithms of length 2 in terms of classical polylogarithms. 
\begin{eqnarray*}
\mathrm{Li}_{1,1}(-1,-1)&=&\frac{1}{2}((\mathrm{Li}_1(-1))^2-\mathrm{Li}_2(1)),\\
\mathrm{Li}_{2,1}(1,-1) &=& -\frac{1}{4}(2\mathrm{Li}_2(1)\mathrm{Li}_1(-1) +\mathrm{Li}_3(1)),\\
\mathrm{Li}_{1,2}(-1,1)&=&\frac{1}{2} (3\mathrm{Li}_2(1)\mathrm{Li}_1(-1) +2\mathrm{Li}_3(1)),\\
\mathrm{Li}_{1,1}(1,-1)&=&\frac{(\mathrm{Li}_1(-1))^2}{2},\\
\mathrm{Li}_{2,1}(-1,-1) &=& \frac{1}{8} (8\mathrm{Li}_2(1)\mathrm{Li}_1(-1) +5\mathrm{Li}_3(1)),\\
\mathrm{Li}_{1,2}(1,1)&=& \mathrm{Li}_3(1),\\
\mathrm{Li}_{1,2}(-1,-1)&=& \frac{1}{8} (-12\mathrm{Li}_2(1)\mathrm{Li}_1(-1) -13\mathrm{Li}_3(1)),\\
\mathrm{Li}_{1,2}(1,-1)&=&\frac{\mathrm{Li}_3(1)}{8}.
\end{eqnarray*}
Applying the previous identities to the expression for $I$ gives
\begin{eqnarray*}
I &=&\frac{1}{2}(\mathrm{Li}_1(-1))^3-\frac{1}{2}\mathrm{Li}_2(1)\mathrm{Li}_1(-1)+\frac{9}{4}\mathrm{Li}_2(1)\mathrm{Li}_1(-1) +\frac{9}{8}\mathrm{Li}_3(1)\\
&&-\frac{3}{4}\mathrm{Li}_2(1)\mathrm{Li}_1(-1) -\frac{1}{2}\mathrm{Li}_3(1)+\frac{(\mathrm{Li}_1(-1))^3}{2}\\
&&- 5\mathrm{Li}_2(1)\mathrm{Li}_1(-1) -\frac{25}{8}\mathrm{Li}_3(1) -\mathrm{Li}_{3}(1)\\
&&+\frac{1}{3}(\mathrm{Li}_1(-1))^3+6\mathrm{Li}_2(1)\mathrm{Li}_1(-1) +\frac{13}{2}\mathrm{Li}_3(1)-\frac{1}{2}\mathrm{Li}_{3}(1)\\
&=&\frac{4}{3}(\mathrm{Li}_1(-1))^3+2\mathrm{Li}_2(1)\mathrm{Li}_1(-1) +\frac{5}{2}\mathrm{Li}_3(1).\\
\end{eqnarray*}
We may now write the expression in terms of values of the zeta function and logarithms.
\[I=-\frac{4}{3}\log^3 2-2\zeta(2)\log 2 +\frac{5}{2}\zeta(3).\]
This shows the required identity for the formula (5).
\qed

The previous Theorem may be completed with the trivial statement 
\[m(x+y+2) = \log 2.\]

In fact, the motivation for setting $c=2$ is that this is the precise point where the family of polynomials $x+y+c$ reaches the unit torus singularly. In classical Mahler measure, those polynomials are among the simplest to compute the Mahler measure, and the same is true in higher Mahler measures.

\subsection{A family related with Dyson integrals}
Consider the following family of polynomials
\begin{eqnarray*}
P_N(x_1,\dots,x_N)&=& \prod_{1 \le h \neq j \le N} \left( 1- \frac{x_h}{x_j} \right) = \prod_{h<j} \left( 2-\frac{x_h}{x_j} - \frac{x_j}{x_h} \right) \\
&=& 2^{N(N-1)} \prod_{h<j} \sin^2 \pi(\theta_h - \theta_j),
\qquad (x_h = \e^{2 \pi \ii \theta_h}).
\end{eqnarray*}
Then we have the following result due to Dyson:
\begin{eqnarray*}
Z(k,P_N)&=& \int_0^1 \cdots \int_0^1 P_N(\e^{2\pi \ii \theta_1}, \dots, \e^{2\pi \ii \theta_N})^k
\dd\theta_1 \cdots \dd \theta_N 
= \frac{(Nk)!}{(k!)^N}.
\end{eqnarray*}

Incorporating this identity into the formula for the zeta Mahler measure we obtain
\begin{eqnarray*}
Z(s,1+\lambda P_N)
&=& \int_0^1 \cdots \int_0^1 (1+\lambda P_N)^s d\theta_1 \cdots d \theta_N \\
&=& \sum_{k=0}^\infty \binom{s}{k} Z(k,P_N) \lambda^k = \sum_{k=0}^\infty \binom{s}{k} \frac{(Nk)!}{(k!)^N} \lambda^k \\
&=& {}_N F_{N-1}
\left(\left. \begin{array}{c}
-s, \frac{1}{N}, \frac2N,\dots,\frac{N-1}{N} \\ 1,\dots,1
\end{array}
\right| \frac{\lambda}{N^N}
\right). \\
\end{eqnarray*}
As always, we may use the expression of zeta to compute higher Mahler measures. By Lemma \ref{property} {\em(ii)},
\begin{eqnarray*}
m(1+\lambda P_N)
&=& \sum_{k=1}^\infty \frac{(-1)^{k-1}}{k} Z(k,P_N) \lambda^k = \sum_{k=1}^\infty \frac{(-1)^{k-1}}{k} \frac{(Nk)!}{(k!)^N} \lambda^k, \\
m_2(1+\lambda P_N)
&=& \sum_{k=1}^\infty \frac{(-1)^{k}}{k} \left(1+\cdots+\frac{1}{k-1} \right) Z(k,P_N) \lambda^k \\
&=& \sum_{k=1}^\infty \frac{(-1)^{k}}{k} \left(1+\cdots+\frac{1}{k-1} \right) 
\frac{(Nk)!}{(k!)^N} \lambda^k.
\end{eqnarray*}
In particular, for $N=2$,
\begin{eqnarray*}
m(1+\lambda P_2) &=& \sum_{k=1}^\infty \frac{(-1)^{k-1}}{k} \binom{2k}{k} \lambda^k, \\
m_2(1+\lambda P_2) &=& \sum_{k=2}^\infty \frac{(-1)^k}{k} \left(1+\cdots+\frac{1}{k-1} \right) 
\binom{2k}{k} \lambda^k.
\end{eqnarray*}
These correspond to the higher Mahler measures of $1+ \lambda(x+x^{-1}+y+y^{-1})$.

\medskip
\noindent{\bf Acknowledgements}: We  would like to thank Fernando Rodriguez-Villegas for helpful discussions.



\begin{flushleft}
Nobushige Kurokawa\\ Department of Mathematics, Tokyo Institute of Technology
2-12-1 Oh-Okayama, Meguro, Tokyo, 152-8551, Japan\\ 
\texttt{kurokawa@math.titech.ac.jp}
\end{flushleft}

\begin{flushleft}
Matilde Lal\'{\i}n\\ Department of Mathematical and Statistical Sciences, University of Alberta,
Edmonton, AB T6G 2G1, Canada\\
\texttt{mlalin@math.ualberta.ca}
\end{flushleft}

\begin{flushleft}Hiroyuki Ochiai\\
Department of Mathematics, Nagoya University
Furo, Chikusa, Nagoya 464-8602, Japan \\
\texttt{ochiai@math.nagoya-u.ac.jp}
\end{flushleft}

\end{document}